\documentclass[11pt]{article}
\usepackage{amsfonts,amsmath,amssymb,amscd,amsthm,mathrsfs}
\title{RIGIDITY OF SUBMANIFOLDS WITH PARALLEL MEAN CURVATURE IN SPACE FORMS\footnote{
\indent Research supported by the NSFC, Grant No.
11071211; the Trans-Century Training Programme
\newline \indent Foundation for Talents by the Ministry of
Education of China.}}
\author{HONG-WEI XU AND JUAN-RU GU}
\date{}
\begin{document}

\maketitle
\begin{abstract}
Let $M$ be an $n(\geq3)$-dimensional oriented compact submanifold
with parallel mean curvature in the simply connected space form
$F^{n+p}(c)$ with $c+H^2>0$, where $H$ is the mean curvature of $M$.
We prove that if the Ricci curvature of $M$ satisfies $Ric_{M}\geq(n-2)(c+H^2),$ then $M$ is either a
totally umbilic sphere, the Clifford hypersurface
$S^{m}\big(\frac{1}{\sqrt{2(c+H^2)}}\big)\times
S^{m}\big(\frac{1}{\sqrt{2(c+H^2)}}\big)$ in
$S^{n+1}(\frac{1}{\sqrt{c+H^2}})$ with $n=2m$, or
$\mathbb{C}P^{2}(\frac{4}{3}(c+H^2))$ in
$S^7(\frac{1}{\sqrt{c+H^2}})$. In particular, if
$Ric_{M}>(n-2)(c+H^2),$ then $M$ is a totally umbilic sphere.
\end{abstract}
 \section{Introduction}
 \hspace*{6mm}The investigation of rigidity of submanifolds with parallel mean curvature attracts a lot of attention
 of differential geometers. After the
pioneering work on compact minimal submanifolds in a sphere due to
Simons \cite{Simons}, Lawson \cite{Lawson} and Chern-do
Carmo-Kobayashi \cite{Chern} obtained a classification of
$n$-dimensional oriented compact minimal submanifolds in $S^{n+p}$
whose squared norm of the second fundamental form satisfies $S\le
n/(2-1/p)$.

It was partially extended to submanifolds with parallel mean
curvature in a sphere by Okumura \cite{Okumura,Okumura1}, Yau
\cite{Yau} and others. In 1990, Xu \cite{Xu0} proved the generalized
Simons-Lawson-Chern-do Carmo-Kobayashi theorem for compact
submanifolds with parallel mean curvature in a sphere.\\\\
\textbf{Theorem A.} \emph{Let $M$ be an $n$-dimensional oriented
compact submanifold with parallel mean curvature in an
$(n+p)$-dimensional unit sphere $S^{n+p}$. Denote by $S$ and $H$ the
squared norm of the second fundamental form and the mean curvature
of $M$. If $S\leq C(n,p,H),$ then $M$ is either a totally umbilic
sphere, a Clifford hypersurface in $S^{n+1}(r)$, or the Veronese
surface in $S^{4}(\frac{1} {\sqrt{1+H^{2}}})$. Here the constant
$C(n,p,H)$ is defined by
$$C(n,p,H)=\left\{\begin{array}{llll} \alpha(n,H),&\mbox{\ $for$\ } p=1, \mbox{\ $or$\ } p=2 \mbox{\  $and$\ }
 H\neq0,\\
 \frac{n}{2-\frac{1}{p}}, &\mbox{\ $for$\ } p\geq2 \mbox{\ $and$\ } H=0, \\
 \min\Big\{\alpha(n,H),\frac{n+nH^2}{2-\frac{1}{p-1}}+nH^2\Big\},&\mbox{\ $fo$r\ } p\geq 3 \mbox{\  $and$\ }
 H\neq0,
 \end{array} \right.$$
where$$\alpha(n,H)=n+ \frac{n^{3}}{2(n-1)}H^{2} -
\frac{n(n-2)}{2(n-1)}\sqrt{n^{2}H^{4}+4(n-1)H^{2}}.$$}

In 1991, A. M. Li and J. M. Li \cite{Li} improved Simons' pinching
constant for $n$-dimensional compact minimal submanifolds in
$S^{n+p}$ to $\max\{\frac{n}{2-1/p},\frac{2}{3}n\}$. Using Li-Li's
matrix inequality \cite{Li}, Xu \cite{Xu} improved the pinching
constant $C(n,p,H)$ in Theorem A to
$$C'(n,p,H)=\left\{\begin{array}{llll} \alpha(n,H),&\mbox{\ for\ } p=1, \mbox{\ or\ } p=2 \mbox{\ and\ }
 H\neq0,\\
\min\Big\{\alpha(n,H),\frac{1}{3}(2n+5nH^2)\Big\},&\mbox{\
otherwise.\ }
\end{array} \right.$$

The rigidity theorem for compact minimal submanifolds with pinched
sectional curvature in a sphere was initiated by Yau \cite{Yau},
then by Itoh \cite{Itoh2}, and finally by Gu and Xu \cite{GX}. It
was extended to compact submanifolds with parallel mean curvature in
space forms by Shen, Han and the authors \cite{GX,Shen,Xu2}.

In 1979, Ejiri \cite{Ejiri} obtained the following rigidity theorem
for $n(\geq 4)$-dimensional oriented compact simply connected
minimal
submanifolds with pinched Ricci curvatures in a sphere.\\\\
\textbf{Theorem  B.} \emph{Let $M$ be an $n(\geq 4)$-dimensional
  oriented compact simply connected minimal submanifold in an
$(n+p)$-dimensional unit sphere $S^{n+p}$. If the Ricci curvature of
$M$ satisfies $Ric_M \geq n-2,$ then $M$ is either the totally
geodesic submanifold $S^n$, the Clifford torus
$S^{m}\big(\sqrt{\frac{1}{2}}\big)\times
S^{m}\big(\sqrt{\frac{1}{2}}\big)$ in $S^{n+1}$ with $n=2m$, or
$\mathbb{C}P^{2}(4/3)$ in $S^7$. Here $\mathbb{C}P^{2}(4/3)$ denotes
the $2$-dimensional complex projective space minimally immersed into
$S^7$ with constant
holomorphic sectional curvature $\frac{4}{3}$.}\\

The pinching constant above is the best possible in
even dimensional cases. It's better than the pinching constants of
Simons \cite{Simons} and Li-Li \cite{Li} in the sense of the average of Ricci curvatures. The following problem
seems very attractive, which has been open for many years.\\\\
\textbf{Problem A.} \emph{Is it possible to generalize Ejiri's
rigidity theorem for minimal submanifolds to the cases of
submanifolds with parallel mean curvature in a
sphere?}\\\\
\hspace*{5mm}In 1987, Sun \cite{Sun} gave a partial answer to the
problem above and showed that if $M$ is an $n(\geq 4)$-dimensional
compact oriented submanifold with parallel mean curvature in
$S^{n+p}$ and its Ricci curvature is not less than
$\frac{n(n-2)}{n-1}(1+H^2)$, then $M$ is a totally umbilic sphere.
Afterward, Shen \cite{ShYB} and Li \cite{Li2} extended Ejiri's
rigidity theorem to the case of 3-dimensional compact minimal
submanifolds in a sphere.

Let $F^{n+p}(c)$ be an $(n+p)$-dimensional simply connected space
form with constant curvature c. Recently, Xu and Tian \cite{XT}
proved a vanishing theorem for the fundamental group of a
submanifold, which says that if $M$ is an $n(\geq 3)$-dimensional
    compact submanifold in the space form
  $F^{n+p}(c)$ with $c\ge0$, and if the Ricci curvature of
$M$ satisfies $Ric_M >\frac{n-1}{2}c+\frac{n^2H^2}{8},$ then $M$ is
simply connected. This implies that the assumption that $M$ is
simply connected in the Ejiri rigidity theorem can be taken off.
Therefore, we have the following refined version of the Ejiri rigidity theorem. \\\\
\textbf{Theorem C.} \emph{Let $M$ be an $n(\geq 3)$-dimensional
  oriented compact minimal submanifold in an
$(n+p)$-dimensional unit sphere $S^{n+p}$. If the Ricci curvature of
$M$ satisfies $Ric_M \geq n-2,$ then $M$ is either the totally
geodesic submanifold $S^n$, the Clifford torus
$S^{m}\big(\sqrt{\frac{1}{2}}\big)\times
S^{m}\big(\sqrt{\frac{1}{2}}\big)$ in $S^{n+1}$ with $n=2m$, or
$\mathbb{C}P^{2}(4/3)$ in $S^7$. Here $\mathbb{C}P^{2}(4/3)$ denotes
the $2$-dimensional complex projective space minimally immersed into
$S^7$ with constant
holomorphic sectional curvature $\frac{4}{3}$.}\\

Motivated by Theorem C, we would like to propose the following problem.\\\\
\textbf{Problem B.} \emph{Is it possible to generalize the refined
version of the Ejiri rigidity theorem for minimal submanifolds to
the cases of submanifolds with parallel mean curvature in space forms?}\\

The purpose of the present paper is to give affirmative
answers to Problems A and B. More precisely, we will prove
 the following rigidity theorem for submanifolds with parallel mean curvature in space forms.\\\\
\textbf{Main Theorem.} \emph{Let $M$ be an $n(\geq3)$-dimensional
 oriented compact submanifold with parallel mean curvature in the
space form $F^{n+p}(c)$ with  $c+H^2>0$.
 If  $$Ric_{M}\geq(n-2)(c+H^2),$$
then M is either a totally umbilic sphere
$S^n(\frac{1}{\sqrt{c+H^2}})$, a Clifford hypersurface
$S^{m}\big(\frac{1}{\sqrt{2(c+H^2)}}\big)\times
S^{m}\big(\frac{1}{\sqrt{2(c+H^2)}}\big)$ in the totally umbilic
sphere $S^{n+1}(\frac{1}{\sqrt{c+H^2}})$ with $n=2m$, or
$\mathbb{C}P^{2}(\frac{4}{3}(c+H^2))$ in
$S^7(\frac{1}{\sqrt{c+H^2}})$. Here
$\mathbb{C}P^{2}(\frac{4}{3}(c+H^2))$ denotes the $2$-dimensional
complex projective space minimally immersed into
$S^7(\frac{1}{\sqrt{c+H^2}})$ with constant
holomorphic sectional curvature $\frac{4}{3}(c+H^2)$.}\\\\

\section{Notation and lemmas}\hspace*{6mm}Throughout this paper,
let $M^{n}$ be an $n$-dimensional compact Riemannian manifold
isometrically immersed in an $(n+p)$-dimensional complete and
simply connected space form $F^{n+p}(c)$ with constant curvature $c$. We
shall make use of the following convention on the range of indices:
$$ 1\leq A,B,C,\ldots\leq n+p;\ 1\leq i,j,k,\ldots\leq n;\ n+1\leq
\alpha,\beta,\gamma,\ldots\leq n+p.$$ Choose a local field of
orthonormal frames \{$e_{A}$\} in $F^{n+p}(c)$ such that, restricted
to $M$, the $e_{i}$'s are tangent to \emph{M}. Let \{$\omega _{A}$\}
and \{$\omega _{AB}$\} be the dual frame field and the connection
1-forms of $F^{n+p}(c)$ respectively. Restricting these forms to
\emph{M}, we have
\begin{eqnarray}&&\omega_{\alpha i}=\sum_{j} h^{\alpha}_{ij}\omega_{j},  \, \,
h^{\alpha}_{ij}=h^{\alpha}_{ji},\nonumber\\
&&h=\sum_{\alpha,i,j}h^{\alpha}_{ij}\omega_{i}\otimes\omega_{j}\otimes
e_{\alpha},\,\,\xi=\frac{1}{n}\sum_{\alpha,i}h^{\alpha}_{ii}e_{\alpha},\nonumber
\\
&&R_{ijkl}=c(\delta_{ik}\delta_{jl}-\delta_{il}\delta_{jk})+\sum_{\alpha}(h^{\alpha}_{ik}h^{\alpha}_{jl}-h^{\alpha}_{il}h^{\alpha}_{jk}), \\
&&R_{\alpha\beta kl}=\sum_{i}(h^{\alpha}_{ik}h^{\beta}_{il}-h^{\alpha}_{il}h^{\beta}_{ik}),
\end{eqnarray}
where $h$, $\xi$, $R_{ijkl}$ and $R_{\alpha\beta kl}$ are the second fundamental form, the mean
curvature vector, the curvature tensor and the normal curvature tensor
of $M$.
Denote by $Ric(u)$ the Ricci curvature of $M$ in direction of $u\in
UM$. From the Gauss equation, we have
\begin{equation}Ric(e_i)=(n-1)c
+\sum_{\alpha}[h_{ii}^{\alpha}h_{jj}^{\alpha}
                          -(h_{ij}^{\alpha})^2].\end{equation}
For a matrix $A=(a_{ij})$, we denote by $N(A)$ the square of the
norm of $A$, i.e.,
$$N(A)=tr(AA^{T})=\sum a^{2}_{ij}.$$
We define
$$S=|h|^{2}, \ H=|\xi|, \ H_{\alpha}=(h^{\alpha}_{ij})_{n\times
n}.$$ Then the scalar curvature $R$ of $M$ is given by
\begin{equation}R=n(n-1)c+n^{2}H^{2}-S.\end{equation}
Let $M$ be a submanifold with parallel
mean curvature vector $\xi$. Choose $e_{n+1}$ such that it is parallel to $\xi$, and
\begin{equation}
trH_{n+1}= nH, \,\, trH_\alpha= 0, \mbox{\ for\ } \alpha \neq n+1.\end{equation}
Set
$$
 S_H=
trH^2_{n+1},\ \ S_I= \sum_{\alpha\neq n+1} trH^2_\alpha.
$$
The following lemma will be used in the proof of our results.\\\\
\textbf{Lemma 1(\cite{Yau}).} \emph{If $M^n$ is a submanifold with parallel mean curvature in $F^{n+p}(c)$, then
either $H\equiv0$, or H is non-zero constant and $H_{n+1}H_\alpha =H_\alpha
H_{n+1}$ for all $\alpha$.}\\\\
  We denote the first and the second covariant derivatives of
$h^{\alpha}_{ij}$ by $h^{\alpha}_{ijk}$ and $h^{\alpha}_{ijkl}$
respectively.  The Laplacian $\Delta h^{\alpha}_{ij}$ of $h^{\alpha}_{ij}$ is defined by $\Delta h^{\alpha}_{ij}=\sum_{k}
h^{\alpha}_{ijkk}$. Following \cite{Yau}, we have
\begin{equation}
\Delta h^{n+1}_{ij}=
\sum_{k,m}(h^{n+1}_{mk}R_{mijk}+h^{n+1}_{im}R_{mkjk}).\end{equation}\\

\section{Proof of Main Theorem}
\hspace*{6mm}To verify Main Theorem, we need to prove the following theorem.\\\\
 \textbf{Theorem 1.} \emph{Let $M$ be an $n(\geq3)$-dimensional oriented
compact submanifold with parallel mean curvature $(H\neq0)$ in the
space form $F^{n+p}(c)$.
 If  $$Ric_{M}\geq(n-2)(c+H^2),$$
where $c+H^2>0$, then M is pseudo-umbilical.}\\\\
\textbf{Proof.} By the Gauss equation (1) and (6), we have
\begin{eqnarray}
\frac{1}{2}\Delta
S_{H}&=&\sum_{i,j,k}(h^{n+1}_{ijk})^{2}+\sum_{i,j}h^{n+1}_{ij}\Delta
h^{n+1}_{ij}\nonumber\\
&=&\sum_{i,j,k}(h^{n+1}_{ijk})^{2}+\sum_{i,j,k,m}h^{n+1}_{ij}h^{n+1}_{km}\Big[(\delta_{mj}\delta_{ik}-
\delta_{mk}\delta_{ij})c+\sum_{\alpha}(h^{\alpha}_{mj}h^{\alpha}_{ik}-h^{\alpha}_{mk}\nonumber\\
&&h^{\alpha}_{ij})\Big]+\sum_{i,j,k,m}h^{n+1}_{ij}h^{n+1}_{im}\Big[(\delta_{mj}\delta_{kk}-
\delta_{mk}\delta_{jk})c+\sum_{\alpha}(h^{\alpha}_{mj}h^{\alpha}_{kk}-h^{\alpha}_{mk}h^{\alpha}_{jk})\Big]\nonumber\\
&=&\sum_{i,j,k}(h^{n+1}_{ijk})^{2}+nc\sum_{i,j}(h^{n+1}_{ij})^{2}-\Big[\sum_{i,j}(h^{n+1}_{ij})^{2}\Big]^{2}-n^{2}cH^{2}\nonumber\\
&&+nH\sum_{i,j,k}h^{n+1}_{ij}h^{n+1}_{jk}h^{n+1}_{ki}-\sum_{\alpha\neq
n+1}\Big[\sum_{i,j}(h^{n+1}_{ij}-H\delta_{ij})h^{\alpha}_{ij}\Big]^{2}.\end{eqnarray}Let $\{e_i\}$ be a frame diagonalizing the matrix $H_{n+1}$ such that $h^{n+1}_{ij}=\lambda_i^{n+1}\delta_{ij}$, for all $i,j$.
Set
\begin{eqnarray*}&&f_k=\sum_{i}(\lambda_i^{n+1})^k,\\
&&\mu_i^{n+1}=H-\lambda_i^{n+1},\,\, i=1,2,...,n,\\
&&B_k=\sum_i(\mu_i^{n+1})^k.
\end{eqnarray*}
Then
\begin{eqnarray*}
&&B_1=0,\,\, B_2=S_H-nH^2,\\
&&B_3=3HS_H-2nH^3-f_k.
\end{eqnarray*}
This together with (7) implies that
\begin{eqnarray}
\frac{1}{2}\Delta{S_H}&=&\sum_{i,j,k}(h^{n+1}_{ijk})^2+ncS_H-S_H^2-n^2cH^2+nHf_3-\sum_{\alpha\neq n+1}\Big(\sum_{i}\mu_i^{n+1}h^\alpha_{ii}\Big)^2\nonumber\\
&=&\sum_{i,j,k}(h^{n+1}_{ijk})^2+ncS_H-S_H^2-n^2cH^2\nonumber\\
&&+nH(3HS_H-2nH^3-B_3)-\sum_{\alpha\neq n+1}\Big(\sum_{i}\mu_i^{n+1}h^\alpha_{ii}\Big)^2\nonumber\\
&=&\sum_{i,j,k}(h^{n+1}_{ijk})^2+B_2[nc+2nH^2-S_H]-nHB_3-\sum_{\alpha\neq n+1}\Big(\sum_{i}\mu_i^{n+1}h^\alpha_{ii}\Big)^2.
\end{eqnarray}
Let $d$ be the infimum of the Ricci curvature of $M$. Then we have
\begin{equation}
Ric(e_i)=(n-1)c+nH\lambda_i^{n+1}-(\lambda_i^{n+1})^2-\sum_{\alpha\neq
n+1,j}(h_{ij}^{\alpha})^2\geq d.
\end{equation}
This implies that
 \begin{equation}
S-nH^2\leq n[(n-1)(c+H^2)-d],
\end{equation}
and
\begin{equation}
(n-2)H(\lambda_i^{n+1}-H)-(\lambda_i^{n+1}-H)^2+(n-1)(c+H^2)-\sum_{\alpha\neq n+1,j}(h_{ij}^{\alpha})^2-d\geq0.
\end{equation}
It follows from (11) that $$
H(\lambda_i^{n+1}-H)\geq\frac{(\lambda_i^{n+1}-H)^2}{n-2}+\frac{\sum_{\alpha\neq
n+1,j}(h_{ij}^{\alpha})^2}{n-2}+\frac{d}{n-2}-\frac{n-1}{n-2}(c+H^2).$$
So,
\begin{eqnarray}
-nHB_3&\geq&\frac{n}{n-2}\sum_i(\mu_i^{n+1})^4+\frac{n}{n-2}
\sum_{\alpha\neq
n+1}\sum_{i,j}(h_{ij}^{\alpha})^2(\mu_i^{n+1})^2\nonumber\\
&&+\frac{n}{n-2}[d-(n-1)(c+H^2)]B_2. \end{eqnarray} From (8) and
(12), we get
\begin{eqnarray}
\frac{1}{2}\Delta{S_H}&\geq&\sum_{i,j,k}(h^{n+1}_{ijk})^2+B_2\Big\{nc+2nH^2-S_H+\frac{n}{n-2}[d-(n-1)(c+H^2)]\Big\}\nonumber\\
&&+\frac{n}{n-2}\sum_i(\mu_i^{n+1})^4+\sum_{\alpha\neq n+1}\Big[\frac{n}{n-2}\sum_{i}(h_{ii}^{\alpha})^2(\mu_i^{n+1})^2-\Big(\sum_{i}\mu_i^{n+1}h^\alpha_{ii}\Big)^2\Big]\nonumber\end{eqnarray}
\begin{eqnarray}
&\geq&\sum_{i,j,k}(h^{n+1}_{ijk})^2+B_2\Big\{nc+2nH^2-S_H+\frac{n}{n-2}[d-(n-1)(c+H^2)]\Big\}\nonumber\\
&&+\frac{B_2^2}{n-2}-\frac{n-3}{n-2}\sum_{\alpha\neq n+1}\Big(\sum_{i}\mu_i^{n+1}h^\alpha_{ii}\Big)^2\nonumber\\
 &\geq&\sum_{i,j,k}(h^{n+1}_{ijk})^2+B_2\Big\{nc+nH^2-\frac{n-3}{n-2}(S-nH^2)\nonumber\\
 &&+\frac{n}{n-2}[d-(n-1)(c+H^2)]\Big\}. \end{eqnarray}
This together with (10) implies that
\begin{eqnarray}
\frac{1}{2}\Delta{S_H}&\geq&\sum_{i,j,k}(h^{n+1}_{ijk})^2+\frac{n}{n-2}B_2\{(n-2)(c+H^2)\nonumber\\
 &&-(n-3)[(n-1)(c+H^2)-d]+[d-(n-1)(c+H^2)]\}\nonumber\\
 &=&\sum_{i,j,k}(h^{n+1}_{ijk})^2+nB_2[d-(n-2)(c+H^2)].
 \end{eqnarray}
 By the assumption, we have $d\ge(n-2)(c+H^2).$ This together with (14) and the maximum principal implies that $S_H$ is a
 constant, and
\begin{equation}(S_H-nH^2)[d-(n-2)(c+H^2)]=0.\end{equation}
Suppose that $S_H\neq nH^2$. Then $d=(n-2)(c+H^2)$. We consider the
following two cases:

  (i) If $n=3$, then the inequalities in (13) and (14) become equalities.
  Thus, we have
  $$h^{\alpha}_{ij}=0, \,\,\,  \mbox{for\,\,}  \alpha\neq n+1, \,\,\,  i\neq j,$$
\begin{equation}|\mu^{n+1}_i|=|\mu^{n+1}_j|, \,\,\,  \mu^{n+1}_i=h^{\alpha}_{ii},
\,\,\,  \mbox{for\,\,} \alpha\neq n+1, \,\,\, 1\le i,j\le
n.\end{equation} This implies $\mu^{n+1}_i=0, \,i=1,2,\cdots ,n$. It
follows from Gauss
  equation that $c+H^2=0.$ This contradicts with
assumption.

  (ii) If $n\geq4$, then the inequalities in (13) and (14) become equalities and we have
  $$Ric_M\equiv(n-2)(c+H^2),$$
 $$h^{\alpha}_{ij}=0, \,\,\,  \mbox{for\,\,}  \alpha\neq n+1, \,\,\,  i\neq j,$$
\begin{equation}|\mu^{n+1}_i|=|\mu^{n+1}_j|, \,\,\,  \mu^{n+1}_i=h^{\alpha}_{ii},  \,\,\,
\mbox{for\,\,} \alpha\neq n+1, \,\,\, 1\le i,j\le n.\end{equation}
It follows from Gauss equation that $\mu^{n+1}_i=0$ and $c+H^2=0.$
This contradicts with assumption.

Therefore, $S_H=nH^2$, i.e., $M$ is a pseudo-umbilical submanifold. This completes the proof of Theorem 1.\\

The following result due to Yau \cite{Yau} will be used in the proof of our main theorem.\\\\
\textbf{Theorem 2.} \emph{Let $N^{n+p}$ be a conformally flat
manifold. Let $N_1$ be a subbundle of the normal bundle of $M^n$
with fiber dimension k. Suppose $M$ is umbilical with respect to
$N_1$ and $N_1$ is parallel in the normal bundle. Then M lies in an
$(n+p-k)$-dimensional umbilical submanifold $N'$ of N such that the
fiber of $N_1$ is everywhere perpendicular to $N'$.} \\\\
\textbf{ Proof of Main Theorem.} When $H=0$, the assertion follows from Theorems C.

When $p=1$ and $H\neq0$, we get the conclusion from Theorem 1.

Now we assume that $p\geq 2$ and $H\neq 0$. It follows from the
assumption and Theorem 1 that $M$ is pseudo-umbilical. It is seen
from Theorem 2 that $M$ lies in an $(n+p-1)$-dimensional totally
umbilic submanifold $F^{n+p-1}(\tilde{c})$ of $F^{n+p}(c)$, i.e.,
the isometric immersion from $M$ into ${F}^{n+p}(c)$ is given by
$$i\circ\varphi: M \rightarrow F^{n+p-1}(\tilde{c})\rightarrow
 F^{n+p}(c),$$ where $\varphi:M^{n}\rightarrow F^{n+p-1}(\tilde{c})$ is
a isometric immersion with  mean curvature vector $\xi_1$, and $i:
F^{n+p-1}(\tilde{c})\rightarrow
 F^{n+p}(c)$
is the totally umbilic submanifold with mean curvature vector
$\xi_2$. Denote by $h_2$ the second fundamental form of isometric
immersion $i$. Set \begin{equation}H_1=|\xi_1|, \,\,\,\,
H_2=|\xi_2|.\end{equation} We know that $\xi=\xi_1+\eta$, where
$\eta=\frac{1}{n}\sum_{i}h_2(e_i,e_i)$ and \{$e_{i}$\} is a local
orthonormal frame field in $M$. Since $\xi_1\perp\xi$, and
$\eta\parallel\xi$, we obtain $\xi_1=0$, and $\eta=\xi$. Noting that
$F^{n+p-1}(\tilde{c})$ is a totally umbilic submanifold in
$F^{n+p}(c)$, we have $|\eta|=H_2$. Thus,
\begin{equation}H^2=H_1^2+|\eta|^2=H_2^2.\end{equation} This together with the Gauss equation
implies that
\begin{equation}\tilde{c}=c+H^2.\end{equation}
Hence, $M$ is an oriented compact minimal submanifold in
$S^{n+p-1}(\frac{1}{\sqrt{c+H^2}})$. It follows from  Theorem C that
$M$ is either a totally umbilic sphere
$S^n(\frac{1}{\sqrt{c+H^2}})$, a Clifford hypersurface
$S^{m}\big(\frac{1}{\sqrt{2(c+H^2)}}\big)\times
S^{m}\big(\frac{1}{\sqrt{2(c+H^2)}}\big)$ in the totally umbilic
sphere $S^{n+1}(\frac{1}{\sqrt{c+H^2}})$ with $n=2m$, or
$\mathbb{C}P^{2}(\frac{4}{3}(c+H^2))$ in
$S^7(\frac{1}{\sqrt{c+H^2}})$. This completes the proof of Main
Theorem.\\\\
\hspace*{5mm}As a consequence of Main Theorem, we get the following\\\\
\textbf{Corollary 1.} \emph{Let $M^n$ be an $n(\geq3)$-dimensional
oriented compact submanifold with parallel mean curvature in the
space form $F^{n+p}(c)$ with  $c+H^2>0$.
 If  $$Ric_{M}>(n-2)(c+H^2),$$
 then M is the totally umbilic sphere $S^n(\frac{1}{\sqrt{c+H^2}})$.}\\\\
\hspace*{5mm}Motivated our main theorem, we would like to propose
the following differentiable rigidity theorem for submanifolds in space forms.\\\\
\textbf{Conjecture A.} \emph{Let $M$ be an $n(\geq3)$-dimensional
compact oriented submanifold in the space form $F^{n+p}(c)$ with $c+H^2>0$.
 If  $$Ric_{M}\geq(n-2)(c+H^2),$$
 then M is diffeomorphic to either the standard
$n$-sphere $S^n$, the Clifford hypersurface
$S^{m}\big(\frac{1}{\sqrt{2}}\big)\times
S^{m}\big(\frac{1}{\sqrt{2}}\big)$ in $S^{n+1}$ with $n=2m$, or
$\mathbb{C}P^{2}$. In particular, if  $Ric_{M}>(n-2)(c+H^2)$, then $M$ is diffeomorphic to $S^n$.}\\\\
\hspace*{5mm}To get an affirmative answer to Conjecture A, we hope
to prove the
following conjecture on the Ricci flow.\\\\
\textbf{Conjecture B.} \emph{Let $(M,g_0)$ be an
$n(\geq4)$-dimensional compact submanifold in an $(n+p)$-dimensional
space form $F^{n+p}(c)$ with $c+H^2>0$. If the Ricci curvature of
$M$ satisfies
$$Ric_{M} >
(n-2)(c+H^2),$$
then the normalized Ricci flow with initial metric
$g_0$
$$\frac{\partial}{\partial t}g(t) = -2Ric_{g(t)} +\frac{2}{n}
r_{g(t)}g(t),$$ exists for all time and converges to a constant
curvature metric as $t\rightarrow\infty$. Moreover, $M$ is
diffeomorphic to a spherical space form. In particular, if M is
simply connected, then M is diffeomorphic to $S^n$.}\\

Making use of the convergence results of Hamilton \cite{Hamilton}
and Brendle \cite{Brendle} for Ricci flow and the nonexistence
theorem of stable currents due to Lawson-Simons \cite{Lawson2} and
Xin \cite{Xin}, Xu and Tian \cite{XT} gave the partial
affirmative answers to Conjectures A and B.\\
\hspace*{5mm}Recently, Andrews and Baker \cite{Andrews}, Liu, Xu, Ye
and Zhao \cite{LXYZ} obtained the convergence theorems for the mean
curvature flow of higher codimension in Euclidean spaces. Motivated
by our main theorem and the conjectures in \cite{GX,LXYZ}, we would
like to propose the following conjecture on the mean
curvature flow in higher codimensions.\\\\
\textbf{Conjecture C.} \emph{Let $F_0:M\rightarrow F^{n+p}(c)$ be an
$n$-dimensional compact submanifold in an $(n+p)$-dimensional space
form $F^{n+p}(c)$ with $c+H^2>0$. If the Ricci curvature of $M$
satisfies
$$Ric_{M} >
(n-2)(c+H^2),$$ then the mean curvature flow
\begin{eqnarray*}
\label{MCF}\left\{
\begin{array}{ll}
\frac{\partial}{\partial t}F(x,t)=n\xi(x,t), \,\, x\in M, \, t\ge0,  \\
F(\cdot,0)=F_0(\cdot),
\end{array}\right.
\end{eqnarray*}
has a unique smooth solution $F : M \times [0, T) \rightarrow
F^{n+p}(c)$ on a finite maximal time interval, and $F_t(\cdot)$
converges uniformly to a round point $q\in F^{n+p}(c)$ as $t
\rightarrow T$. In particular, $M$ is diffeomorphic to $S^n$.}\\\\

Hong-Wei Xu

Center of Mathematical Sciences\

Zhejiang University\

Hangzhou 310027\

China

E-mail address: xuhw@cms.zju.edu.cn\\\\

Juan-Ru Gu

Center of Mathematical Sciences\

Zhejiang University\

Hangzhou 310027\

China

E-mail address: gujr@cms.zju.edu.cn

\end{document}